\documentclass[twoside]{article}
\usepackage{a4,txmac,%
mypic,%
nocaphead,%
amssymb,times,mathptm}

\makeatletter

\advance\oddsidemargin by -1.7cm
\advance\evensidemargin by -1.7cm
\advance\textwidth by 3.4cm

\def\mynewtheo#1#2{%
\newtheorem{@#1}{#2}[section]%
\newenvironment{#1}{\begin{@#1}\rm}{\end{@#1}}}
 
\mynewtheo{lemma}{Lemma}
\mynewtheo{exer}{Exercise}
\mynewtheo{theo}{Theorem}
\mynewtheo{rem}{Remark}
\mynewtheo{defi}{Definition}
\mynewtheo{conj}{Conjecture}
\mynewtheo{corr}{Corollary}
\mynewtheo{prop}{Proposition}
\mynewtheo{question}{Question}
\mynewtheo{exam}{Example}

\def\inx{\mathop {\operator@font ind}\mathord{}}
\def\mpb{\mathop {\operator@font mpb}\mathord{}}
\def\mwf{\mathop {\operator@font mwf}\mathord{}}
\def\spn{\mathop {\operator@font span}\mathord{}}
\def\len{\mathop {\operator@font len}\mathord{ }}

\newenvironment{theorem}{\begin{theo}}{\end{theo}}

\def\eqref#1{\mbox{(\protect\reference{#1}})}
\def\proof{\@ifnextchar[{\@proof}{\@proof[\unskip]}}
\def\@proof[#1]{\noindent{\bf Proof #1.}\enspace}

\begin{document}

\author{A. Stoimenow\footnotemark[1]\\[2mm]
\small Research Institute for Mathematical Sciences, \\
\small Kyoto University, Kyoto 606-8502, Japan\\
\small e-mail: {\tt stoimeno@kurims.kyoto-u.ac.jp}\\
\small WWW: {\hbox{\web|http://www.kurims.kyoto-u.ac.jp/~stoimeno/|}}
}

{\def\thefootnote{}
\footnotetext[1]{This is a preprint. I would be grateful
  for any comments and corrections. Current
  version: \today\ \ \ First version: \makedate{13}{4}{2007}}
\def\thefootnote{\fnsymbol{footnote}}
\footnotetext[1]{Financial support by the 21st Century COE Program
is acknowledged.}
}

\title{\large\bf
\uppercase{Tait's conjectures and}\\[2mm] 
\uppercase{odd crossing number amphicheiral knots}
}

\date{}

\maketitle

\renewcommand{\section}{%
   \@startsection
         {section}{1}{\z@}{1.9ex \@plus -1ex \@minus -.2ex}%
               {.9ex \@plus.2ex}{\large\bf}%
}

\renewcommand{\subsubsection}{%
   \@startsection
         {subsubsection}{1}{\z@}{-1.5ex \@plus -1ex \@minus -.2ex}%
               {1ex \@plus.2ex}{\large\bf}%
}
\renewcommand{\@seccntformat}[1]{\csname the#1\endcsname .
\quad}

\def\chrd#1#2{\picline{1 #1 polar}{1 #2 polar}}
\def\arrow#1#2{\picvecline{1 #1 polar}{1 #2 polar}}

\def\labch#1#2#3{\chrd{#1}{#2}\picputtext{1.3 #2 polar}{$#3$}}
\def\labar#1#2#3{\arrow{#1}{#2}\picputtext{1.3 #2 polar}{$#3$}}
\def\labbr#1#2#3{\arrow{#1}{#2}\picputtext{1.3 #1 polar}{$#3$}}

\def\rottab#1#2{%
\expandafter\advance\csname c@table\endcsname by -1\relax
\centerline{%
\rbox{\centerline{\vbox{\setbox1=\hbox{#1}%
\centerline{\mbox{\hbox to \wd1{\hfill\mbox{\vbox{{%
\caption{#2}}}}\hfill}}}%
\vskip9mm
\centerline{
\mbox{\copy1}}}}%
}%
}%
}

\def\GD{\szCD{6mm}}
\def\szCD#1#2{{\let\@nomath\@gobble\small\diag{#1}{2.4}{2.4}{
  \picveclength{0.27}\picvecwidth{0.1}
  \pictranslate{1.2 1.2}{
    \piccircle{0 0}{1}{}
    #2
}}}}
\def\CD{\szCD{4mm}}

\def\point#1{{\picfillgraycol{0}\picfilledcircle{#1}{0.08}{}}}
\def\labpt#1#2#3{\pictranslate{#1}{\point{0 0}\picputtext{#2}{$#3$}}}
\def\vrt#1{{\picfillgraycol{0}\picfilledcircle{#1}{0.09}{}}}

\def\@dcont{}
\def\svCD#1{\ea\glet\csname #1\endcsname\@dcont}
\def\rsCD#1{\ea\glet\ea\@dcont\csname #1\endcsname\ea\glet
\csname #1\endcsname\relax}

\def\addCD#1{\ea\gdef\ea\@dcont\ea{\@dcont #1}}
\def\drawCD#1{\szCD{#1}{\@dcont}}

\def\noloop{{\diag{0.8cm}{0.5}{1}{\picline{0.25 0}{0.25 1}}}}

\def\vrt#1{{\picfillgraycol{0}\picfilledcircle{#1}{0.09}{}}}

\def\ReidI#1#2{
  \diag{0.8cm}{0.9}{1}{
    \pictranslate{0.4 0.5}{
      \picscale{#11 #21}{
        \picmultigraphics[S]{2}{1 -1}{
	  \picmulticurve{-6 1 -1 0}{0.5 -0.5}{0.5 0}{0.1 0.3}{-0.2 0.3}
	} 
	\piccirclearc{-0.2 0}{0.3}{90 270}
      }
    }
  }
}
\def\Pos#1#2{{\diag{#1}{1}{1}{#2
\picmultiline{-8 1 -1 0}{0 1}{1 0}
\picmultiline{-8 1 -1 0}{0 0}{1 1}
}}}
\def\Neg#1#2{{\diag{#1}{1}{1}{#2
\picmultiline{-8 1 -1 0}{0 0}{1 1}
\picmultiline{-8 1 -1 0}{0 1}{1 0}
}}}
\def\Nul#1#2{{\diag{#1}{1}{1}{#2
\piccirclearc{0.5 1.4}{0.7}{-135 -45}
\piccirclearc{0.5 -0.4}{0.7}{45 135}
}}}
\def\Inf#1#2{{\diag{#1}{1}{1}{#2
\piccirclearc{0.5 1.4 x}{0.7}{135 -135}
\piccirclearc{0.5 -0.4 x}{0.7}{-45 45}
}}}
\def\pos{\Pos{0.5em}{\piclinewidth{10}}}
\def\neg{\Neg{0.5em}{\piclinewidth{10}}}
\def\nul{\Nul{0.5em}{\piclinewidth{10}}}
\def\iinf{\Inf{0.5em}{\piclinewidth{10}}}

\def\pt#1{{\picfillgraycol{0}\picfilledcircle{#1}{0.1}{}}}
\def\ppt#1{{\picfillgraycol{0}\picfilledcircle{#1}{0.06}{}}}

\def\@curvepath#1#2#3{%
  \@ifempty{#2}{\piccurveto{#1 }{@stc}{@std}#3}%
    {\piccurveto{#1 }{#2 }{#2  #3  0.5 conv}
    \@curvepath{#3}}%
}
\def\curvepath#1#2#3{%
  \piccurve{#1 }{#2 }{#2 }{#2  #3  0.5 conv}%
  \picPSgraphics{/@stc [ #1  #2  -1 conv ] $ D /@std [ #1  ] $ D }%
  \@curvepath{#3}%
}

\def\@opencurvepath#1#2#3{%
  \@ifempty{#3}{\piccurveto{#1 }{#1 }{#2 }}%
    {\piccurveto{#1 }{#2 }{#2  #3  0.5 conv}\@opencurvepath{#3}}%
}
\def\opencurvepath#1#2#3{%
  \piccurve{#1 }{#2 }{#2 }{#2  #3  0.5 conv}%
  \@opencurvepath{#3}%
}

\def\@br#1#2#3#4#5#6{
  \pictranslate{#1}{
     \picrotate{#2}{
        #6{
          \picmultigraphics[S]{2}{-1 1}{
            \picmulticurve{-7 1 -1 0}
	      {#3 -.5 * #3 .5 *}{#3 -.2 * #3 .5 *}
              {#3 .2 * #3 -.5 *}{#3 .5 * #3 -.5 *}
          }
        }
        \picmultigraphics[S]{2}{-1 1}{
	   \@ifempty{#5}{
             \picmultigraphics[S]{2}{1 -1}{
               \picline{#3 .5 * d}{#4 .5 * #3 .5 *}
	     }
	   }{
	     \picscale{1 -1}{\picline{#3 .5 * d}{#4 .5 * #3 .5 *}}
	     \picvecline{#3 .5 * d}{#4 #3 d + + 6 : #3 .5 *}
	     \piclineto{#4 .5 * #3 .5 *}
	   }
        }
     }
  }
}

\def\lbr#1#2#3#4{\@br{#1}{#2}{#3}{#4}{}{\picscale{-1 1}}}
\def\lvbr#1#2#3#4{\@br{#1}{#2}{#3}{#4}{a}{\picscale{-1 1}}}
\def\rbr#1#2#3#4{\@br{#1}{#2}{#3}{#4}{}{}}
\def\rvbr#1#2#3#4{\@br{#1}{#2}{#3}{#4}{a}{}}

\long\def\@makecaption#1#2{%
   \vskip 10pt
   {\let\label\@gobble
   \let\ignorespaces\@empty
   \xdef\@tempt{#2}%
   }%
   \ea\@ifempty\ea{\@tempt}{%
   \setbox\@tempboxa\hbox{%
      \fignr#1#2}%
      }{%
   \setbox\@tempboxa\hbox{%
      {\fignr#1:}\capt\ #2}%
      }%
   \ifdim \wd\@tempboxa >\captionwidth {%
      \rightskip=\@captionmargin\leftskip=\@captionmargin
      \unhbox\@tempboxa\par}%
   \else
      \hbox to\captionwidth{\hfil\box\@tempboxa\hfil}%
   \fi}%
\def\fignr{\small\sffamily\bfseries}%
\def\capt{\small\sffamily}%

\newdimen\@captionmargin\@captionmargin2\parindent\relax
\newdimen\captionwidth\captionwidth\hsize\relax

\def\nin{\not\in}
\def\bQ{{\Bbb Q}}
\def\bC{{\Bbb C}}
\def\bR{{\Bbb R}}
\def\bN{{\Bbb N}}
\def\bZ{{\Bbb Z}}
\def\cE{{\cal E}}
\def\cK{{\cal K}}
\def\cI{{\cal I}}
\def\cR{{\cal R}}
\def\cV{{\cal V}}
\def\cO{{\cal O}}
\def\hD{{\hat D}}
\def\cD{{\cal D}}
\def\cX{{\cal X}}
\def\hK{{\hat K}}
\def\bm{\bar t'_2}
\let\bt\bm 

\def\pr{\text{\rm pr}\,}
\def\ncap{\not\mathrel{\cap}}
\def\|{\mathrel{\kern1.5pt\Vert\kern1.5pt}}
\def\lra{\longrightarrow}
\def\llra{\longleftrightarrow}
\let\ds\displaystyle
\let\reference\ref
\def\so{\Rightarrow}
\let\x\exists
\let\fa\forall
\let\es\enspace

\let\x\exists
\let\sg\sigma
\let\tl\tilde
\let\ap\alpha
\let\dl\delta
\let\Dl\Delta
\let\be\beta
\let\gm\gamma
\let\Gm\Gamma
\let\nb\nabla
\let\Lm\Lambda
\let\sm\setminus
\let\vn\varnothing
\let\dt\det
\def\tG{\tl G}
\def\tD{\tl D}
\def\tS{\tl S}

\let\eps\varepsilon
\let\ul\underline
\let\ol\overline
\def\md{\min\deg}
\def\Md{\max\deg}
\def\Mcf{\max{\operator@font cf}}
\let\Mc\Mcf
\def\sgn{{\operator@font sgn}}

\def\ssim{\stackrel{\ds \sim}{\vbox{\vskip-0.2em\hbox{$\scriptstyle *$}}}}

\def\bysame{\same[\kern2cm]\,}
\def\qed{\hfill\@mt{\Box}}
\def\@mt#1{\ifmmode#1\else$#1$\fi}

\def\proof{\@ifnextchar[{\@proof}{\@proof[\unskip]}}
\def\@proof[#1]{\noindent{\bf Proof #1.}\enspace}

\def\myfrac#1#2{\raisebox{0.2em}{\small$#1$}\!/\!\raisebox{-0.2em}{\small$#2$}}
\newcommand{\mybr}[2]{\text{$\Bigl\lfloor\mbox{%
\small$\displaystyle\frac{#1}{#2}$}\Bigr\rfloor$}}
\def\mybrtwo#1{\mbox{\mybr{#1}{2}}}

\def\epsfs#1#2{{\ifautoepsf\unitxsize#1\relax\else
\epsfxsize#1\relax\fi\epsffile{#2.eps}}}
\def\epsfsv#1#2{{\vcbox{\epsfs{#1}{#2}}}}
\def\vcbox#1{\setbox\@tempboxa=\hbox{#1}\parbox{\wd\@tempboxa}{\box
  \@tempboxa}}

\def\@test#1#2#3#4{%
  \let\@tempa\go@
  \@tempdima#1\relax\@tempdimb#3\@tempdima\relax\@tempdima#4\unitxsize\relax
  \ifdim \@tempdimb>\z@\relax
    \ifdim \@tempdimb<#2%
      \def\@tempa{\@test{#1}{#2}}%
    \fi
  \fi
  \@tempa
}

\def\go@#1\@end{}
\newdimen\unitxsize
\newif\ifautoepsf\autoepsftrue

\unitxsize4cm\relax
\def\epsfsize#1#2{\epsfxsize\relax\ifautoepsf
  {\@test{#1}{#2}{0.1 }{4   }
		{0.2 }{3   }
		{0.3 }{2   }
		{0.4 }{1.7 }
		{0.5 }{1.5 }
		{0.6 }{1.4 }
		{0.7 }{1.3 }
		{0.8 }{1.2 }
		{0.9 }{1.1 }
		{1.1 }{1.  }
		{1.2 }{0.9 }
		{1.4 }{0.8 }
		{1.6 }{0.75}
		{2.  }{0.7 }
		{2.25}{0.6 }
		{3   }{0.55}
		{5   }{0.5 }
		{10  }{0.33}
		{-1  }{0.25}\@end
		\ea}\ea\epsfxsize\the\@tempdima\relax
		\fi
		}

\let\diagram\diag

\def\boxed#1{\diagram{1em}{1}{1}{\picbox{0.5 0.5}{1.0 1.0}{#1}}}

\def\rato#1{\hbox to #1{\rightarrowfill}}
\def\arrowname#1{{\enspace
\setbox7=\hbox{F}\setbox6=\hbox{%
\setbox0=\hbox{\footnotesize $#1$}\setbox1=\hbox{$\to$}%
\dimen@\wd0\advance\dimen@ by 0.66\wd1\relax
$\stackrel{\rato{\dimen@}}{\copy0}$}%
\ifdim\ht6>\ht7\dimen@\ht7\advance\dimen@ by -\ht6\else
\dimen@\z@\fi\raise\dimen@\box6\enspace}}

\def\contr{\diagram{1em}{0.6}{1}{\piclinewidth{35}%
\picstroke{\picline{0.5 1}{0.2 0.4}%
\piclineto{0.6 0.6}\picveclineto{0.3 0}}}}

\newcounter{pp}%
\newenvironment{mylist}[1]{%
\begin{list}{#1{pp})}%
{\usecounter{pp}\setlength{\labelwidth}{4mm}%
\setlength{\leftmargin}{0.6cm}\setlength{\itemsep}{1mm}%
\setlength{\topsep}{1mm}}%
\gdef\myitem{\item\xdef\@currentlabel{#1{pp})}}%
}{\end{list}}

\def\abstractname{}

\parskip5pt plus 1pt minus 2pt
\parindent\z@

{\let\@noitemerr\relax
\vskip-2.7em\kern0pt\begin{abstract}
\noindent{\bf Abstract.}\enspace
We give a brief historical overview of the Tait conjectures,
made 120 years ago in the course of his pioneering work in 
tabulating the simplest knots, and solved a century later
using the Jones polynomial. We announce the solution, again
based on a substantial study of the Jones polynomial, 
of one (possibly his last remaining?) 
problem of Tait, with the construction of amphicheiral knots 
of almost all odd crossing numbers. An application to the
non-triviality problem for the Jones polynomial is also
outlined.
\\[1mm]
(secondary)
\end{abstract}
}
\vspace{1mm}


\section{The first knot tables}

Knot theory took its origins in the late 19th century.
At that time, W.~Thomson (``Lord Kelvin''), P.~G.~Tait
and J.~Maxwell propagated the vortex-atom theory, in 
an attempt to explain the structure of the universe. 
They believed that a super-substance, ether, makes up
all of matter, and atoms are knotted tubes of ether.
Knotting is hereby understood as tying a piece of rope,
and then identifying its both ends so that the tying cannot
be any more undone.

Thus, in the realm of constructing a periodic table
of elements, Tait began the catalogisation of the
``simplest'' \em{knots}. He depicted knots (as we still 
do today) by \em{diagrams}, consisting of a (smooth) 
plane curve with transverse self-intersections, or 
\em{crossings}. At each crossing one of the two strands 
passes over the other. The above notion of simplicity refers 
to the number of crossings of the diagram. Tait's list aimed 
at presenting, among (knots with) diagrams of few
crossings, each knot by exactly one diagram. Alternatively 
we can define the \em{crossing number} of a knot as the 
minimal crossing number of all its diagrams, and say that 
we seek the list of knots with given (small) crossing 
number. We also like knots represented by different
diagrams in the list to be inequivalent, in the sense that
one cannot turn a (closed) piece of rope knotted the one 
way into one knotted the other way, without cutting the rope.
The simplest knots are shown in figure \ref{fkn}. The
leftmost one, of crossing number 0, is the \em{trivial} knot or 
\em{unknot}. It has some special importance, much like the unit
element in a group.

\begin{figure}[htb]
\[
\begin{array}{c@{\quad}c@{\kern-5mm}c@{}c}
\diag{6mm}{3}{3}{
  \piclinewidth{18}
  \piccircle{1.5 1.5}{1.5}{}
} & \diag{1cm}{4}{3}{ \trefoil } & \diag{1cm}{4}{3}{ \rtrefoil }
  & \diag{1cm}{4}{3}{ \fig8knot } 
\\
\mbox{unknot} & \mbox{left-hand trefoil} & \mbox{right-hand trefoil} &
\mbox{figure-8 knot}
\end{array}
\]
\caption{\label{fkn}}
\end{figure}

Tait completed the list up to 7 crossings. Little, Kirkman,
later Conway \cite{Conway} and others took over and continued 
his work. In the modern computer age, tables have reached
the knots of 17 crossings, with millions of entries, even
though Tait's vortex-atom theory has long been dismissed. An
account on knot tabulation is given, with emphasis on its history,
in \cite{Silver}, and from a more contemporary point of view
in \cite{HTW,Hoste}.

\section{Tait's conjectures}

Tait's accomplishment allows us to call him with some right the
first knot theorist. Yet Tait worked mainly by intuition. He
had at his time no rigorous way of showing knots inequivalent.
Tait's reasoning is not easy to interpret precisely nowadays,
had it been formulated in a language quite different from (and far
less developed than) our present. Nonetheless he evidently observed
several phenomena, which, apart from knot tabulation, would become
a legacy to his successors.

It seems to remain unclear whether Tait was convinced certain
properties to hold for all, or just for alternating knots. A 
knot is \em{alternating}, if in some (alternating) diagram the 
curve passes crossings interchangingly over-under like
$\diag{7mm}{1.5}{1}{
  \picmultiline{-6 1 -1 0}{0.55 0 x}{0.55 1 x}
  \picmultiline{-6 1 -1 0}{0.4 0.1}{0.4 1}
  \picmultiline{-6 1 -1 0}{1.1 0.1}{1.1 1}
  \picmultiline{-6 1 -1 0}{0.55 1.5 x}{0.55 1 x}
}
$, i.e. not containing 
$\,\diag{7mm}{1.5}{1}{
  \picmultiline{-6 1 -1 0}{0.55 0 x}{0.55 1.5 x}
  \picmultiline{-6 1 -1 0}{0.4 0.1}{0.4 1}
  \picmultiline{-6 1 -1 0}{1.1 0.1}{1.1 1}
}\,
$ or
$\,\diag{7mm}{1.5}{1}{
  \picmultiline{-6 1 -1 0}{0.4 0.1}{0.4 1}
  \picmultiline{-6 1 -1 0}{1.1 0.1}{1.1 1}
  \picmultiline{-6 1 -1 0}{0.55 0 x}{0.55 1.5 x}
}\,
$. The knots
in figure \ref{fkn} are such. In fact, this is true for all
knots up to 7 crossings, catalogued by Tait, and at least
for a large portion of the slightly more complicated ones he 
was shown by his successors in his lifetime (even though it is
known now that alternation is a rare property for generic crossing
numbers \cite{Thistle3}). Thus, certainly Tait was guided by
evidence from alternating diagrams. Their occurrence in the 
tables suggested to him

\begin{conj}(Tait's conjecture I)
A \em{reduced}, i.e.~not of the form\ \ $\diag{5mm}{4}{2}{
  \picPSgraphics{/default@brdrm [-8.5 1 -1 0] D}
  \picrotate{-90}{\rbraid{-1 2}{1 1.4}}
  \picfillgraycol{0.4 0.5 0.8}
  \picscale{1 -1}{
    \picfilledcircle{0.7 -1}{0.8}{}
  }
  \picfillgraycol{0.7 0.2 0.6}
  \picfilledcircle{3.3 1}{0.8}{}
}$\,\,, alternating diagram has minimal crossing number
(for the knot it represents).
\end{conj}

For the next problem, we need to define the \em{writhe}.
If one equips (the curve of) a knot diagram with an orientation, 
then each crossing looks, if observed from an appropriate angle,
locally like $\diag{5mm}{1}{1}{
\picmultivecline{0.28 1 -1.0 0}{1 0}{0 1}
\picmultivecline{0.28 1 -1.0 0}{0 0}{1 1}
}$ (\em{positive} crossing) or $\diag{5mm}{1}{1}{
\picmultivecline{0.28 1 -1.0 0}{0 0}{1 1}
\picmultivecline{0.28 1 -1.0 0}{1 0}{0 1}
}$ (\em{negative} 
crossing). The writhe is the difference between the number of 
former crossings and the number of latter. (An easy observation
shows that the writhe is the same for either orientation.)

\begin{conj}(Tait's conjecture II)
Minimal crossing number diagrams of the same (alternating?)
knot have the same writhe.
\end{conj}

For alternating diagrams, he conjectured more precisely the
following:

\begin{conj}(Tait's conjecture III)
Alternating  diagrams of the same knot are related by
a sequence of \em{flypes}:
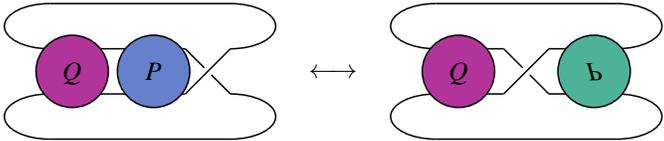
\[
\diag{6mm}{6}{3}{
    \pictranslate{3 1.5}{
       \picmultigraphics[S]{2}{1 -1}{
           \picmultiline{0.22 1 -1.0 0}{2 -0.5}{1 0.5}
           \picmultigraphics[S]{2}{-1 1}{
                \picellipsearc{-2 -1.0}{1 0.5}{90 270}
           }
           \picline{-2 -1.5}{2 -1.5}
           \picline{-2 -0.5}{1 -0.5}
      }
      \picfillgraycol{0.4 0.5 0.8}
      \picfilledcircle{0.3 0}{0.8}{$P$}
      \picfillgraycol{0.7 0.2 0.6}
      \picfilledcircle{-1.5 0}{0.8}{$Q$}
   }
}\quad\llra\quad
\diag{6mm}{6}{3}{
    \pictranslate{3 1.5}{
       \picmultigraphics[S]{2}{1 -1}{
           \picmultiline{0.22 1 -1.0 0}{0.5 -0.5}{-0.5 0.5}
           \picmultigraphics[S]{2}{-1 1}{
                \picellipsearc{-2 -1.0}{1 0.5}{90 270}
                \picline{-2 -0.5}{-0.5 -0.5}
           }
           \picline{-2 -1.5}{2 -1.5}
      }
      \picscale{1 -1}{
           \picfillgraycol{0.3 0.7 0.6}
           \picfilledcircle{1.5 0}{0.8}{$P$}
      }
      \picfillgraycol{0.7 0.2 0.6}
      \picfilledcircle{-1.5 0}{0.8}{$Q$}
   }
}
\]
\end{conj}

One can easily observe that a flype preserves the writhe,
and so conjecture III implies conjecture II (for alternating
knots\footnote{up to some technical issues of primeness and\label{fn1}
whether all minimal crossing diagrams are alternating. These
issues are settled, but we like to skip them here for simplicity.}).

For Tait's last problem, we consider an \em{amphicheiral} knot.
Such a knot can be turned into its mirror image. From the knots
in figure \ref{fkn}, the unknot is obviously amphicheiral. So 
is the figure-eight knot, as shows a simple exercise. In contrast,
the trefoil is not amphicheiral. In other words, the left-hand
trefoil and its mirror image, the right-hand trefoil, are
two distinct knots (a fact that stubbed
knot theorists for a while, and was first proved by Max Dehn).
Tait was wondering what crossing numbers amphicheiral knots can 
have. The evidence he had in mind can probably be formulated so:

\begin{conj}(Tait's conjecture IV)
Amphicheiral (alternating?) knots have even crossing number.
\end{conj}

Note that (for alternating knots, and with the remark in footnote
\reference{fn1}) this is a consequence of conjecture II, for mirroring
a diagram interchanges positive and negative crossings, and so
negates the writhe.

\section{Reidemeister moves and invariants}

A few decades after their genesis, Tait's knot lists 
would be \em{proved} right. With the work of Alexander, 
Reidemeister and others knot theory began to be put on a
mathematical fundament. Reidemeister showed that three
types of local moves (i.e. moves altering only a fragment
of the diagrams, as shown below) suffice to interrelate
all diagrams of the same knot. 
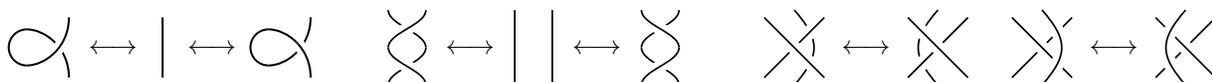
\[
\begin{array}{c@{\kern1cm}c@{\kern1cm}c}
\ReidI{ }{ }\es\llra\,\noloop\,\llra\es\ReidI{ }{-} &
\diag{0.5cm}{1}{2}{
  \piclinewidth{22}
  \picPSgraphics{/default@brdrm [-7.5 1 -1 0] D}
  \rbraid{0.5 0.5}{1 1}
  \lbraid{0.5 1.5}{1 1}
}\es\llra\es
\diag{0.5cm}{1}{2}{
  \piclinewidth{22}
  \picline{0 0}{0 2}\picline{1 0}{1 2}
}\es\llra\es
\diag{0.5cm}{1}{2}{
  \piclinewidth{22}
  \picPSgraphics{/default@brdrm [-7.5 1 -1 0] D}
  \lbraid{0.5 0.5}{1 1}
  \rbraid{0.5 1.5}{1 1}
} &
\diag{0.5cm}{2}{2}{
  \piclinewidth{22}
  \picmulticurve{-8 1 -1 0}{1 0}{1.7 0.8}{1.7 1.2}{1 2}
  \picmultiline{-8 1 -1 0}{0.2 0.2}{1.8 d}
  \picmultiline{-8 1 -1 0}{0.2 1.8}{1.8 0.2}
}\,\llra\,
\diag{0.5cm}{2}{2}{
  \piclinewidth{22}
  \picmulticurve{-8 1 -1 0}{1 0}{0.3 0.8}{0.3 1.2}{1 2}
  \picmultiline{-8 1 -1 0}{0.2 0.2}{1.8 d}
  \picmultiline{-8 1 -1 0}{0.2 1.8}{1.8 0.2}
}\quad
\diag{0.5cm}{2}{2}{
  \piclinewidth{22}
  \picmultiline{-8 1 -1 0}{0.2 0.2}{1.8 d}
  \picmultiline{-8 1 -1 0}{0.2 1.8}{1.8 0.2}
  \picmulticurve{-8 1 -1 0}{1 0}{1.7 0.8}{1.7 1.2}{1 2}
}\,\llra\,
\diag{0.5cm}{2}{2}{
  \piclinewidth{22}
  \picmultiline{-8 1 -1 0}{0.2 0.2}{1.8 d}
  \picmultiline{-8 1 -1 0}{0.2 1.8}{1.8 0.2}
  \picmulticurve{-8 1 -1 0}{1 0}{0.3 0.8}{0.3 1.2}{1 2}
} 
\end{array}
\]
To distinguish two knots thus translates into the question how 
to prove that two diagrams (of these knots) are not connected 
by a sequence of such moves. This is done with the help
of an \em{invariant}, that is, a map
\[
\{\mbox{\ knot diagrams\ }\}\to\mbox{\ ``something''}\,,
\]
whose value in ``something'' does not change (is invariant) 
when the argument (diagram) is changed by a Reidemeister move. 
The question what ``something'' should be is justified. The 
answer is that it would suffice to be any set of objects
whose distinctness is easy to verify, yet which is
large enough to allow the invariant
to take many different values. For us it will be the ring 
$\bZ[t,t^{-1}]$ of Laurent polynomials in one variable with 
integer coefficients. Clearly one can compare coefficients easier
than wondering about a sequence of Reidemeister moves, which 
may be arbitrarily long and pass over arbitrarily complicated
intermediate diagrams. (There is a fundamental, but out of our
focus here, issue how to, and that one can to some extent, control
these sequences \cite{HL}.) An obvious other desirable feature 
of an invariant is that we could easily evaluate it from a diagram. 

Alexander's merit was to construct precisely such an invariant.
The \em{Alexander polynomial} \cite{Alexander} remained (and
still remains) a main theme in knot theory for decades to
come. Let us observe, though, that we actually already 
came across one other knot invariant, the crossing 
number. It is an invariant directly by definition,
since it was defined on the whole Reidemeister move
equivalence class of diagrams. However, this definition
makes it difficult to evaluate from a diagram~-- in contrast 
to Alexander's polynomial, for which its creator, and 
later many others, gave several simple procedures.

\section{The Jones polynomial}

60 years after Alexander, a new chapter of knot theory
was opened by V.~Jones, with the discovery of a successor
to Alexander's polynomial. The developments the \em{Jones 
polynomial} $V$ \cite{Jones} has sparked in the 20 years
since its appearance are impossible even to be vaguely 
sketched in completeness, and go far beyond both the
competence and expository intention of the author here.
Several more concepts, like links, braids, geometric,
Vassiliev and quantum invariants etc., are left 
out for simplicity and length reasons, for which the 
author likes to apologise at this point. A good (though 
still partial, and now no longer very recent) account 
on these issues was given by Birman \cite{Birman}.

Let us recall, though, one of the first achievements
the Jones polynomial became famous with~-- the solution
of Tait's conjectures. Conjectures I, II and IV, for 
alternating knots, were proved by Kauffman \cite{Kauffman2}, 
Murasugi \cite{Murasugi,Murasugi2} and Thistlethwaite \cite{Thistle2}.
We will need a few more words on Kauffman's proof, since it 
uses a calculation procedure for $V$ called \em{state model}.
This state model also gives a very elementary proof
that the Jones polynomial is an invariant. (Kauffman
had previously developed a similar model for the Alexander 
polynomial, too \cite{Kauffman3}.)

Based on Kauffman's state model, Lickorish and Thistlethwaite 
\cite{LickThis} defined a \em{(semi)adequate} knot and diagram. An
advantage of this concept is that alternating diagrams/knots are
adequate. Many details would better be skipped here, but let us
say that a diagram is adequate/semiadequate if it is $+$adequate 
and/or $-$adequate, and that taking the mirror image of a diagram
transforms the property for $+$ into the one for $-$. A knot
with some of these properties is defined as one that has a diagram
with the feature of the same name. Thistlethwaite extended 
the proof of the Tait's three conjectures to the class of adequate
knots \cite{Thistle}, applying a 2-variable generalisation of
$V$, the Kauffman polynomial \cite{Kauffman}.

For a general (non-alternating) knot, in case of conjecture II,
Tait's intuition had been proved misleading. In the 1970's,
K.~Perko \cite{Perko} observed that there are two 10 crossing
knots in the tables \cite[appendix]{Rolfsen}, which are
equivalent, even though their 10 crossing diagrams have different 
writhe. This duplication had remained unnoticed for quite 
a while, possibly due to the belief in Tait's conjecture.
Remedying this error (and a few other long-remained ones) 
in the tables still causes some confusion in their use.

Tait's conjecture III was settled a few years after
the others by Menasco and Thistlethwaite \cite{MenThis}, 
mostly using geometric techniques, though again with some 
(now subordinate) appearance of the Jones polynomial.

\section{The crossing numbers of amphicheiral knots}

Kauffman, Murasugi and Thistlethwaite's proof of Tait's 
conjecture IV shows that an alternating knot $K$ of odd 
crossing number and its mirror image $!K$ always have
distinct Jones polynomials. (Let us in contrast remark that $K$ 
and $!K$ have the same Alexander polynomial for \em{every} knot $K$.)
In the opposite direction, their work allows to easily find a(n 
alternating) amphicheiral knot of every even crossing number
at least 4. But their results could not decide what crossing 
number non-alternating amphicheiral knots can have. The main
aim of this note is to announce the complete solution to Tait's
(last?) problem.

\begin{theorem}\label{tM}
For each odd natural number $n\ge 15$, there exists an
amphicheiral knot of crossing number $n$.
\end{theorem}

Similarly to Perko's knot, a particular instance disproving Tait's
conjecture IV for non-alternating knots was found accidentally:
Hoste and Thistlethwaite, in the course of routine knot tabulation,
discovered an amphicheiral 15 crossing knot. (Their compilational
work had previously shown that there are no amphicheiral knots
of odd crossing numbers up to 13.) Settling the other crossing
numbers is a major problem, 
though, since exhaustive enumeration is no longer a feasible 
approach~-- we face the above noticed difficulty that we do not 
know (generally) how to determine the crossing number. A few 
other methods are known, but all they fail on such examples. 
Thus the way to our result is rather far, and below we will 
conclude by just giving a brief outline of the proof. The details
will appear in a separate (long) paper. We also mention another 
application of our approach, which addresses the non-triviality 
problem for the Jones polynomial.

\section{Semiadequacy invariants and the 
  non-triviality problem}

The Alexander polynomial was, from its very beginning, connected
to topological features of knots. The situation is rather different
for its successor. The problem to give a topological meaning to
the Jones polynomial has bothered many knot theorists ever
since this invariant appeared. So far we still find 
ourselves in the embarrassing state where we ``can quickly
fill pages with the coefficients and exponents of $V$
for not-too-complicated knots without having the slightest
idea what they mean'' (\cite[end of \S 3]{Birman}). Similarly
unsolved, and intriguing, remains the problem, formulated by
Jones, if his polynomial detects the unknot. Again, ``our lack
of knowledge about this problem is in striking contrast to
the control mathematicians now have over the Alexander
polynomial: understanding its topological meaning, we also
know precisely how to construct knots with Alexander polynomial
$1$'' (ibid., rem. (iii) after theorem 3, \S 8; see, though, also
\cite{EKT}).

In an attempt to gain more insight into the appearance of the
Jones polynomial, the author, and (up to minor interaction,
independently) Dasbach and Lin \cite{DL,DL2}, initiated a
detailed study of some coefficients of $V$ in semiadequate
diagrams. Let us remark here that, while adequacy is only a 
slight extension of alternation, semiadequacy is a rather 
wide extension of adequacy. (For the experts: semiadequate
knots contain completely positive, Montesinos and 3-braid
knots.) Semiadequacy is still a fairly general condition,
yet it helps settle many technical issues. 

For semiadequate knots the first coefficient of the Jones
polynomial is $\pm 1$, almost by definition \cite{LickThis}.
The outcome of our work was that we gained an understanding
of coefficients 2 and 3. Their invariance allows to derive
3 invariant quantities each from a $+$adequate, and similarly
from a $-$adequate diagram, called below \em{semiadequacy
invariants}. Their merit is that they reflect directly
certain features of the diagram, and so we have a precise
idea how a semiadequate diagram with given invariants
must look like. The first of them allows to prove:

\begin{theorem}
Semiadequate knots have non-trivial Jones polynomial.
\end{theorem}

This implies (for experts) the result also for Montesinos 
and 3-braid knots, but it can be proved also for their 
Whitehead doubles, some strongly $n$-trivial knots and 
$k$-almost positive knots with $k\le 3$.

Our three semiadequacy invariants become also, joined
by a relative obtained from the Kauffman polynomial and
Thistlethwaite's results \cite{Thistle}, the main tool for 
the proof of theorem \ref{tM}. For given odd $n\ge 15$, we
start with an amphicheiral knot $K$ that has an $n$ crossing
diagram, which is semiadequate. Luckily, such examples
can be obtained by leaning on Hoste-Thistlethwaite's knot.
The work in \cite{Thistle} shows then that the crossing number
of $K$ is at least $n-1$, and were it $n-1$, a minimal crossing
diagram $D$ would be adequate. Then we have 4 invariants for
both $+$adequacy and $-$adequacy each available. A detailed
study of how an $n-1$ crossing diagram with such invariants 
must look like is necessary to exclude most cases for $D$.
Hereby, among the various generalisations of Thistlethwaite's
knot, one must choose carefully the one whose invariants make
the exclusion argument most convenient (or better to say, feasible
at all). Only a small fraction of possibilities for $D$
remain, which are easy to check, and rule out, by computer.
This allows us to conclude that in fact $D$ cannot exist.


\end{document}